\documentstyle[11pt,amssymb]{article}   
\begin{document}  

\date{} 
\renewcommand{\thefootnote}{}

\def\thebibliography#1{\begin{center}{\normalsize\bf References}\end{center}
 \list{[{\bf \arabic{enumi}}]}{\settowidth\labelwidth{[#1]}
 \leftmargin\labelwidth
 \advance\leftmargin\labelsep
 \usecounter{enumi}}
 \def\newblock{\hskip .11em plus .33em minus .07em}
 \sloppy\clubpenalty4000\widowpenalty4000
 \sfcode`\.=1000\relax}

\title{\vspace*{-1.5cm}
{\large TORUS KNOTS THAT CANNOT BE UNTIED BY TWISTING}}

\author{{\normalsize MOHAMED AIT NOUH}\\
{\small Department of Mathematics, University of Provence}\\[-1mm]
{\small 39 rue F. Joliot-Curie, 13453 Marseille cedex 13, France}\\[-1mm]
{\small e-mail: Mohamed.Aitnouh@cmi.univ-mrs.fr}\\[3mm]
{\normalsize AKIRA YASUHARA }\\
{\small Department of Mathematics, Tokyo Gakugei University}\\[-1mm]
{\small Nukuikita 4-1-1, Koganei, Tokyo 184-8501, Japan}\\
{\small {\em Current address}, October 1, 1999 to September 30, 2001:}\\[-1mm]
{\small Department of Mathematics, The George Washington University}\\[-1mm]
{\small Washington, DC 20052, USA}\\[-1mm]
{\small e-mail: yasuhara@u-gakugei.ac.jp}\\
}

\maketitle

\baselineskip=12pt
\vspace*{-5mm}  
{\small 
\begin{quote}
\begin{center}A{\sc bstract}\end{center} 
We give a necessary condition for a torus knot 
to be untied by a single twisting. 
By using this result, we give infinitely many 
torus knots that cannot be untied by a single twisting. 
\end{quote}}


\baselineskip=16pt
\begin{center}{\bf{1. Introduction}}\end{center}

Throughout this paper, we work in the smooth category. 
All orientable manifolds will be assumed 
to be oriented unless otherwise stated. In particular all knots are oriented. 
For an oriented manifold $M$, $-M$ denotes $M$ with the opposite orientation. 

Let $K$ be a knot in the $3$-sphere $S^3$, 
and $D^2$  a disk intersecting $K$ in 
its interior. Let $\omega=|{\mathrm{lk}}(\partial D^2,L)|$ 
and $n$ an integer. A $-1/n$-Dehn 
surgery along $\partial D^2$ changes $K$ into a new knot $K'$ in $S^3$. 
We say that $K'$ is obtained from $K$ by {\em$(n,\omega)$-twisting} 
(or simply {\em twisting}). 
(The second author calls an $(n,\omega)$-twisting 
a $(-n,\omega)$-twisting in his prior papers \cite{KY}, \cite{MY1} 
and \cite{Yasu}.)
Then we write 
$K'\stackrel{(n,\omega)}{\rightarrow}K$.
Let $\cal T$ denote the set of knots that are obtained from 
a trivial knot by a single twisting.
Y. Ohyama \cite{Ohy} showed that any knot can be untied by two twistings.
This implies that any knot is obtained from a knot in $\cal T$ by 
a single twisting. 

A {\em$(p,q)$-torus knot} $T(p,q)$ is a knot that wraps around 
the standard solid torus in the longitudinal direction $p$ times  
and the meridional direction $q$ times, 
where the linking number of the meridian 
and longitude is equal to 1. Note that $p$ and $q$ are coprime. 
A torus knot $T(p,q)$ $(0<p<q)$ is {\em exceptional} 
if $q\equiv \pm1$ (mod $p$), and 
{\em non-exceptional} if it is not exceptional.

Let $p(\geq 2)$ be an integer. It is not hard to see that 
$T(p,\pm1)\stackrel{(k,p)}{\rightarrow}T(p,kp\pm1)$. 
Since $T(p,\pm1)$ is a trivial knot, $T(p,kp\pm1)$ 
belongs to $\cal T$. 
(In \cite{Mot}, K. Motegi calls $T(p,kp\pm1)$ 
a {\em trivial example} of torus knots that 
belong to $\cal T$.) 
This implies that any exceptional torus knot 
belongs to $\cal T$. 
In particular, all of the knots 
$T(2,q)$, $T(3,q)$, $T(4,q)$ and $T(6,q)$ belong 
to $\cal T$. 
In contrast with this fact, a non-exceptional torus knot 
that belongs to $\cal T$ is not known so far. 
These facts let us hit on the following.

\medskip
{\bf Conjecture.} 
{\em No non-exceptional torus knot belongs to $\cal T$.}

\medskip\noindent 
This conjecture seems likely to be true. 
However a non-exceptional torus knot that is not contained in $\cal T$ 
is not known. So we are faced with the following problem 
before this conjecture. 

\medskip
{\bf Problem.} {\em Is there a torus knot that is not 
contained in $\cal T$?}

\medskip
\noindent
In this paper we give a necessary condition for a non-exceptional 
torus knot to belong to $\cal T$, 
and by using this condition, we give infinitely many non-exceptional 
torus knots that are not contained in $\cal T$. 
K. Miyazaki and the second author \cite{MY1} gave a sufficient condition 
for a knot not to be contained in $\cal T$ and showed that there are 
infinitely many knots that are not contained in $\cal T$.
The sufficient condition given in \cite{MY1} cannot be applied to 
torus knots since it contains the condition that the value of the 
signature is equal to 0. 
(It is known that the signature of a nontrivial torus knot does not 
vanish; see Corollary 2.2 and also see \cite{Prz2} for example.)

For a prime integer $d$, let $\sigma_d(K)$ be the {\em Tristram's 
$d$-signature} of a knot $K$ \cite{Tri}. Note that $\sigma_2(K)$ is the 
same as the signature $\sigma(K)$ in the usual sense \cite{Tro}, 
\cite{Mur}. 

\medskip
{\bf Theorem 1.1.} {\em Let $T(p,q)$ $(0<p<q)$ 
be a non-exceptional torus knot. 
If $T(p,q)$ is obtained from a trivial knot by a single 
$(n,\omega)$-twisting, then 
{\rm (i)} $n=1$, {\rm(ii)} $\omega<q$, 
{\rm (iii)} $\omega>p$ if $\omega$ is even, and 
{\rm(iv)} if $\omega$ is divisible by a prime integer $d$, then 
\[\frac{2[d/2](d-[d/2])}{d^2}\omega^2=-\sigma_d(T(p,q)) 
\mbox{ or } =2-\sigma_d(T(p,q)),\]
where $[x]$ is the greatest integer not exceeding $x$.}

\medskip
{\bf Remark 1.2.}
In \cite{MM}, K. Miyazaki and K. Motegi showed that 
if a non-exceptional torus knot $T(p,q)$ $(0<p<q)$ is obtained from 
a trivial knot by a single $(n,\omega)$-twisting, then $|n|=1$. 
Thus we eliminate the possibility $n=-1$. 

\medskip
By using this theorem, we have the following three results.

\medskip
{\bf Theorem 1.3.}
{\em Let $p$ be an odd integer.
If $p\geq 9$, $p\equiv 1$ or $\equiv 3$ $($mod $8)$, then 
$T(p,p+4)$ does not belong to $\cal T$ }

\medskip
{\bf Remark 1.4.} 
Let $p$ be an odd integer. By the argument similar to 
that in the proof Theorem 1.3, we see that 
if $p\geq 7$, $p\equiv 5$ or $\equiv 7$ (mod 8),
and if $T(p,p+4)$ is obtained from a trivial knot by 
a single $(n,\omega)$-twisting, then 
$n=1$ and $\omega=p+2$. 

\medskip
{\bf Theorem 1.5.}
{\em Let $r$ be an even integer. \\
{\rm(1)} If $r\geq 4$, $2p\geq (r/2+1)^2-r/2$, and  $p\equiv 1$ $($mod $2r)$, 
then $T(p,p+r)$ does not belong to $\cal T$. \\
{\rm(2)} If $r\geq 8$, $2p\geq (r/2-1)^2-r/2$, and $p\equiv -1$ $($mod $2r)$, 
then $T(p,p+r)$ does not belong to $\cal T$.} 

\medskip
{\bf Example 1.6.} Let $r$ be an even integer and $n$ a positive integer. 
By the theorem above, we have the following: 
If $4\leq r\leq 14$ and $p=2nr+1$, or if $8\leq r\leq 20$ and $p=2nr-1$, 
then $T(p,p+r)$ does not belong to $\cal T$. Note that this 
contains the case $p\equiv 1$ (mod 8) of Theorem 1.3.

\medskip
{\bf Remark 1.7.} Let $p(\geq 5)$ be an odd integer.
By the argument similar to that in the proof of 
Theorem 1.5, we see that 
if $T(p,p+2)$ is obtained from a trivial knot by a 
single $(n,\omega)$-twisting, then 
$n=1$ and $\omega=p+1$.

\medskip
Since the knots $T(2,q)$, $T(3,q)$, $T(4,q)$ and $T(6,q)$
are exceptional, $T(5,7)$ is the \lq minimum' non-exceptional 
torus knot, i.e., the {crossing number} of $T(5,7)$ is 
minimum in the crossing numbers of non-exceptional torus knots. 
By Remark 1.7, if $T(5,7)$ is obtained from a trivial knot by a 
single $(n,\omega)$-twisting, then 
$n=1$ and $\omega=6$. 
The authors cannot eliminate the possibility 
$(n,\omega)=(1,6)$. 
So it is still open if $T(5,7)$ belongs to $\cal T$ or not.
Concerning  $T(5,8)$, which is the minimum one except for $T(5,7)$, 
we have the following.

\medskip
{\bf Proposition 1.8.} {$T(5,8)$ does not belong to $\cal T$.}

\bigskip
\begin{center}
{\bf 2. {Signatures of torus knots}}
\end{center}

In this section, we calculate the signatures 
of torus knots.

\medskip
{\bf Proposition 2.1.} {\em Let $T(p,q)$ $(0<p<q)$ be 
a torus knot. Then 
\[\sigma(T(p,q))=
2\sum_{i=1}^{[(p-1)/2]}\left(\left[\frac{(p-2i)q}{2p}\right]-
\left[\frac{(3p-2[p/2]-2i)q}{2p}\right]\right)
+(p-1-2[p/2])(q-1).\]}

\medskip
{\bf Proof.} By \cite[Proposition 1]{Lit}, we have 
$\sigma(T(p,q))=\sigma^{+}-\sigma^-$, 
where
\[\begin{array}{rl}
\sigma^{+}=&\displaystyle
\hspace*{-.5em}\#\left\{(i,j)\in {\Bbb Z}\times {\Bbb Z}
\left|0<i<p,0<j<q, 0<\frac{i}{p}+\frac{j}{q}<\frac{1}{2}\right.\right\}\\
& \displaystyle 
+\#\left\{(i,j)\in {\Bbb Z}\times {\Bbb Z}\left|0<i<p,0<j<q, 
\frac{3}{2}<\frac{i}{p}+\frac{j}{q}<2\right.\right\}
\end{array}\]
and\
\[\sigma^-=\#\left\{(i,j)\in {\Bbb Z}\times {\Bbb Z}
\left|0<i<p,0<j<q, \frac{1}{2}<\frac{i}{p}+\frac{j}{q}<\frac{3}{2}
\right.\right\}.\]
We note that if $p>0$ and $q>0$, then 
\[\begin{array}{rl}
0<i/p+j/q<1/2 &\Leftrightarrow\ -qi/p<j<(p-2i)q/2p,\\ 
3/2<i/p+j/q<2 &\Leftrightarrow\ (3p-2i)q/2p<j<(2p-i)q/p,\\
1/2<i/p+j/q<3/2 &\Leftrightarrow\ (p-2i)q/2p<j<(3p-2i)q/2p,\\
(p-2i)q/2p\leq 0 &\Leftrightarrow\ i\geq p/2, \\
(3p-2i)q/2p\geq q &\Leftrightarrow\ i\leq p/2,
\end{array}\]
and $-qi/p<0$, $(p-2i)q/2p<q$ and $(2p-i)q/p\geq q$ for $0<i<p$.
So we have 
\[\sigma^{+}=\#\left\{j\left|0<i<\frac{p}{2}, 
0<j<\frac{(p-2i)q}{2p}\right.\right\}
+\#\left\{j\left|\frac{p}{2}<i<p, \frac{(3p-2i)q}{2p}<j<q\right.\right\} \]
and
\[\sigma^-=\#\left\{j\left|0<i\leq \frac{p}{2}, \frac{(p-2i)q}{2p}<j<q\right.
\right\}+
\#\left\{j\left|\frac{p}{2}<i<p, 0<j<\frac{(3p-2i)q}{2p}\right.\right\}.\]
Since $0<(p-2i)/p<1$ for $0<i<p/2$, and $p$ and $q$ are coprime, 
$(p-2i)q/2p$ is not an integer for $0<i<p/2$. 
Suppose $(3p-2i)q/2p$ is an integer for some $i$ $(p/2<i<p)$. 
If $q$ is odd, then $(3p-2i)/2p$ is an integer. Since $1/2<(3p-2i)/2p<1$ 
for $p/2<i<p$, this is absurd. Therefore $q$ is even. 
Then $p$ is odd, and hence $i/p$ is an integer. This is a contradicton. 
So $(3p-2i)q/2p$ is not an integer for $p/2<i<p$. 
It follows from that 
\[\sigma^{+}=\#\left\{j\left|0<i<\frac{p}{2}, 0<j\leq
\left[\frac{(p-2i)q}{2p}\right]\right.\right\}
+\#\left\{j\left|\frac{p}{2}<i<p, \left[\frac{(3p-2i)q}{2p}\right]<j<q
\right.\right\}\] 
and 
\[\sigma^-=\#\left\{j\left|0<i\leq \frac{p}{2}, 
\left[\frac{(p-2i)q}{2p}\right]<j<q\right.\right\}+
\#\left\{j\left|\frac{p}{2}<i<p, 0<j\leq
\left[\frac{(3p-2i)q}{2p}\right]\right.\right\}.\]
This implies that
\[
\sigma^{+}=\sum_{i=1}^{[(p-1)/2]}\left[\frac{(p-2i)q}{2p}\right]+
\sum_{i=[p/2]+1}^{p-1}\left(q-1-\left[\frac{(3p-2i)q}{2p}\right]\right) 
\]
and 
\[
\sigma^-=\sum_{i=1}^{[p/2]}\left(q-1-\left[\frac{(p-2i)q}{2p}\right]\right)+
\sum_{i=[p/2]+1}^{p-1}\left[\frac{(3p-2i)q}{2p}\right].\]
Thus it is not hard to see that
\[
\sigma(T(p,q))=
2\sum_{i=1}^{[(p-1)/2]}\left[\frac{(p-2i)q}{2p}\right]-
2\sum_{i=[p/2]+1}^{p-1}\left[\frac{(3p-2i)q}{2p}\right]
+(p-1-2[p/2])(q-1).\]
We note that 
\[\sum_{i=[p/2]+1}^{p-1}\left[\frac{(3p-2i)q}{2p}\right]=
\sum_{i=1}^{[(p-1)/2]}\left[\frac{(3p-2[p/2]-2i)q}{2p}\right].\]
This completes the proof. $\Box$

\medskip
{\bf Corollary 2.2.} {\em Let $T(p,q)$ $(0<p<q)$ be a torus knot. Then 
\[\sigma(T(p,q))\leq -2\left[\frac{p}{2}\right]\left[\frac{q}{2}\right]\].}

{\bf Proof.} 
Suppose $p$ is odd. 
Then, by Propositon 2.1, 
\[\sigma(T(p,q))=
2\sum_{i=1}^{(p-1)/2}\left(\left[\frac{(p-2i)q}{2p}\right]-
\left[\frac{(2p-2i+1)q}{2p}\right]\right).\]
Since
\[\left[\frac{(p-2i)q}{2p}\right]-\left[\frac{(2p-2i+1)q}{2p}\right]
=\left[\frac{(p-2i)q}{2p}\right]-\left[\frac{(p-2i)q}{2p}+\frac{(p+1)q}{2p}
\right]
\leq-\left[\frac{(p+1)q}{2p}\right],\]
we have 
\[\sigma(T(p,q))\leq-(p-1)\left[\frac{(p+1)q}{2p}\right]
\leq-(p-1)\left[\frac{q}{2}\right]=
-2\left[\frac{p}{2}\right]\left[\frac{q}{2}\right].\]
Suppose $p$ is even. Note that $q$ is odd.
Then, by Propositon 2.1, 
\[\sigma(T(p,q))=
2\sum_{i=1}^{p/2-1}\left(\left[\frac{(p-2i)q}{2p}\right]-
\left[\frac{(2p-2i)q}{2p}\right]\right)-q+1.\]
Since
 \[\left[\frac{(p-2i)q}{2p}\right]-\left[\frac{(2p-2i)q}{2p}\right]
=\left[\frac{(p-2i)q}{2p}\right]-
\left[\frac{(p-2i)q}{2p}+\frac{pq}{2p}\right]\leq-\left[\frac{q}{2}\right],\]
we have 
\[\sigma(T(p,q))\leq-2(p/2-1)\left[\frac{q}{2}\right]-q+1=
-\frac{(p-2)(q-1)}{2}-q+1=
-2\left[\frac{p}{2}\right]\left[\frac{q}{2}\right].\]
This completes the proof. $\Box$

\medskip
{\bf Proposition 2.3.} {\em Let $p(>0)$ be an odd integer and 
$r$ $(2\leq r<p)$ an even integer, and $T(p,p+r)$ a torus knot. Then 
\[\sigma(T(p,p+r))=
-\frac{(p-1)(p+r+1)}{2}+2\sum_{i=1}^{r/2}
\left(\left[\frac{(2i-1)p}{2r}\right]-\left[\frac{(2i-1)p+r}{2r}\right]
\right)\]}

\medskip
{\bf Proof.} 
By Proposition 2.1, we have 
\[
\sigma(T(p,p+r))=
2\sum_{i=1}^{(p-1)/2}\left(\left[\frac{(p-2i)(p+r)}{2p}\right]-
\left[\frac{(2p-2i+1)(p+r)}{2p}\right]\right).\]
Note that
\[
\left[\frac{(p-2i)(p+r)}{2p}\right]=\frac{(p-2i-1)}{2}
+\left[\frac{(r+1)}{2}-\frac{ri}{p}\right].\]
Since 
\[
\left[\frac{(r+1)}{2}-\frac{ri}{p}\right]=\left\{
\begin{array}{ll}
\displaystyle \frac{r}{2} & \displaystyle 0<i\leq\frac{p}{2r}, \\
\displaystyle \frac{r}{2}-1 & \displaystyle 
\frac{p}{2r}< i\leq\frac{3p}{2r},\\
\vdots & \ \ \vdots \\
1 & \displaystyle \frac{(r-3)p}{2r}< i\leq\frac{(r-1)p}{2r}, \\
0 & \displaystyle \frac{(r-1)p}{2r}< i\leq \frac{p-1}{2},
\end{array}\right.\]
we have  
\[ 
\sum_{i=1}^{(p-1)/2}\left[\frac{(r+1)}{2}-\frac{ri}{p}\right]=
\left[\frac{p}{2r}\right]+\left[\frac{3p}{2r}\right]+
\cdots+\left[\frac{(r-1)p}{2r}\right].\] 
Meanwhile, we have 
\[
\left[\frac{(2p-2i+1)(p+r)}{2p}\right]=
p-i+\left[r+\frac{1}{2}-\frac{(2i-1)r}{2p}\right].\]
Note that if $p>2r$, then
\[
\left[r+\frac{1}{2}-\frac{(2i-1)r}{2p}\right]=\left\{
\begin{array}{ll}
r & \displaystyle 0<i\leq\frac{p+r}{2r},\\
r-1 & \displaystyle \frac{p+r}{2r}<i\leq\frac{3p+r}{2r},\\
\vdots & \ \ \vdots\\
\displaystyle \frac{r}{2}+1 & 
\displaystyle \frac{(r-3)p+r}{2r}<i\leq\frac{(r-1)p+r}{2r},\\
\displaystyle \frac{r}{2} & 
\displaystyle \frac{(r-1)p+r}{2r}< i\leq \frac{p-1}{2},
\end{array}\right.\]
and if $(r<)p<2r$, then
\[
\left[r+\frac{1}{2}-\frac{(2i-1)r}{2p}\right]=\left\{
\begin{array}{ll}
r & \displaystyle 0<i\leq\frac{p+r}{2r},\\
r-1 & \displaystyle \frac{p+r}{2r}< i\leq\frac{3p+r}{2r},\\
\vdots & \ \ \vdots\\
\displaystyle \frac{r}{2}+2 & 
\displaystyle \frac{(r-5)p+r}{2r}<i\leq\frac{(r-3)p+r}{2r},\\
\displaystyle \frac{r}{2}+1 & 
\displaystyle \frac{(r-3)p+r}{2r}<i\leq \frac{p-1}{2}.
\end{array}\right.\]
This implies 
\[
\sum_{i=1}^{(p-1)/2}\left[r+\frac{1}{2}-\frac{(2i-1)r}{2p}\right]
=\left[\frac{p+r}{2r}\right]+\left[\frac{3p+r}{2r}\right]
+\cdots+\left[\frac{(r-1)p+r}{2r}\right]+\frac{r(p-1)}{4},\] 
if $p>2r$, and  
\[
\sum_{i=1}^{(p-1)/2}\left[r+\frac{1}{2}-\frac{(2i-1)r}{2p}\right]
=\left[\frac{p+r}{2r}\right]+\left[\frac{3p+r}{2r}\right]
+\cdots+\left[\frac{(r-3)p+r}{2r}\right]+\frac{(r+2)(p-1)}{4},\]  
if $p<2r$. If $p<2r$, then since $p$ is odd, 
\[\frac{(r+2)(p-1)}{4}-\frac{r(p-1)}{4}
=\frac{p-1}{2}
=\left[\frac{p+1}{2}-\frac{p}{2r}\right]
=\left[\frac{(r-1)p+r}{2r}\right].\] 
Thus we have the required equation.  $\Box$

\medskip
{\bf Proposition 2.4.} {\em Let $p=2nr\pm1(>0)$ be an integer, 
$r$ $(2\leq r< p)$ an even integer, and $T(p,p+r)$ a torus knot. Then 
\[\sigma(T(p,p+r))=\left\{
\begin{array}{ll}
\displaystyle
-\frac{(p-1)(p+r+1)}{2} & \mbox{ if $p=2nr+1$},\\
\displaystyle
-\frac{(p-1)(p+r+1)}{2}-r & \mbox{ if $p=2nr-1$}.
\end{array}\right.\]}

\medskip
{\bf Proof.} 
We note that
$[(2i-1)p/2r]=(2i-1)n$ and $[((2i-1)p+r)/2r]=(2i-1)n$ for $0<i\leq r/2$ 
if $p=2nr+1$, and 
$[(2i-1)p/2r]=(2i-1)n-1$ and $[((2i-1)p+r)/2r]=(2i-1)n$ for $0<i\leq r/2$ 
if $p=2nr-1$. This and Proposition 2.3 complete the proof. $\Box$

\medskip
{\bf Proposition 2.5.} {\em Let $p(>0)$ be an odd integer and 
$T(p,p+4)$ a torus knot. Then 
\[\sigma(T(p,p+4))=\left\{
\begin{array}{ll}
\displaystyle
-\frac{(p-1)(p+5)}{2} & \mbox{ if $p\equiv 1$, or $\equiv3$ $($mod $8)$},\\
\displaystyle
-\frac{(p-1)(p+5)}{2}-4 & \mbox{ if $p\equiv 5$, or $\equiv7$ $($mod $8)$}.
\end{array}\right.\]}

\medskip
{\bf Proof.} Since 
$\sigma(T(1,5))=0$, $\sigma(T(3,7))=-8$, we may assume that $p>4$. 
By combining Proposition 2.3 and the following, we complete proof.
\[\left[\frac{p}{8}\right]+\left[\frac{3p}{8}\right]
-\left[\frac{p+4}{8}\right]-\left[\frac{3p+4}{8}\right]=
\left\{
\begin{array}{ll}
\displaystyle
0 & \mbox{ if $p\equiv 1$, or $\equiv3$ (mod 8)},\\
\displaystyle
-2 & \mbox{ if $p\equiv 5$, or $\equiv7$ (mod 8).  $\Box$} 
\end{array}\right. \]

\bigskip
\begin{center}
{\bf 3. {Proofs of Theorems 1.1, 1.3, 1.5 and Proposition 1.8}}
\end{center}

Similar results to the following two lemmas, Lemmas 3.1 and 3.2, 
are mentioned in several articles 
\cite{Wei}, \cite{Yam}, \cite{Yasu}, \cite{MY1}, \cite{MY2}, \cite{KY}, etc. 
The first lemma is a spacial case of \cite[Lemma 4.4]{KY}. 
The second one is proven by combining \cite[Example 2]{MY2} and the proof of 
\cite[Lemma 2.3]{MY2}. 

\medskip
{\bf Lemma 3.1.} {\em Let $K_1$ and $K_2$ be knots. Let $M$ be a twice 
punctured $-\varepsilon CP^2$. If 
$K_1\stackrel{(\varepsilon,\omega)}{\rightarrow}K_2$, then there exists 
an annulus $A$ in $M$ such that 
$(\partial M, \partial A)\cong(-S^3,-K_1)\cup(S^3,K_2)$ and $A$ represents 
a homology element $\omega\gamma$, 
where $|\varepsilon|=1$ and 
$\gamma$ is a standard generator 
of $H_2(M,\partial M;{\Bbb Z})$ with 
the intersection number $\gamma\cdot\gamma=-\varepsilon$. $\Box$}

\medskip
{\bf Lemma 3.2.} {\em Let $K_1$ and $K_2$ be knots. Let $M$ be a twice 
punctured $S^2\times S^2$. If $K_1\stackrel{(2n,\omega)}{\rightarrow}K_2$, 
then there exists 
an annulus $A$ in $M$ such that 
$(\partial M, \partial A)\cong(-S^3,-K_1)\cup(S^3,K_2)$ and 
$A$ represents a  homology element $\omega\alpha-n\omega\beta$, where 
$\alpha,\beta$ are standard generators of $H_2(M,\partial M;{\Bbb Z})$ 
with $\alpha\cdot\alpha=\beta\cdot \beta=0$ 
and $\alpha\cdot\beta=1$. $\Box$}

\medskip
The following theorem is originally due to O.Ya. Viro \cite{Vir}. 
It is also obtained by letting $a=[{d}/{2}]$ in the inequality 
of \cite[Remarks(a) on p-371]{Gil} by P. Gilmer. 

\medskip
{\bf Theorem 3.3.} (P.M. Gilmer \cite{Gil}, O.Ya. Viro \cite{Vir})
{\em Let $M$ be a compact, oriented, once punctured $4$-manifold, 
and $K$ a knot in $\partial M$.
Suppose that $K$ bounds a properly embedded, oriented 
surface $F$ in $M$ that represents an element 
$\xi\in H_2(M,\partial M;{\Bbb Z})$. 
If $\xi$ is divisible by a prime integer $d$, then we have
\[\left|\frac{2[d/2](d-[d/2])}{d^2}\xi\cdot\xi
-\sigma(M)-\sigma_d (K)\right|
\leq {\rm dim}H_2(M;{\Bbb Z}_p)+2\,{\rm genus}(F). \ \Box\]}

\medskip
The following is a well known result for $d=2$ \cite{Mur} \cite{Gill}. 
J.H. Przytycki showed it in \cite{Prz}. 
Here we show it by using Theorem 3.3.

\medskip
{\bf Lemma 3.4.} 
{\em Let $K_+$ and $K_-$ be knots. If $K_-$ is obtained from $K_+$ by 
changing a positive crossing into negative one, 
then for any prime integer $d$
\[\sigma_d(K_-)-2\leq\sigma_d(K_+)\leq\sigma_d(K_-).\]}

\medskip
{\bf Proof.} 
It is not hard to see that $K_+\stackrel{(1,0)}{\rightarrow}K_-$. 
This implies that 
$K_+\#(-\overline{K_+})\stackrel{(1,0)}{\rightarrow}
K_-\#(-\overline{K_+})$, where $-\overline{K_+}$ is the refrected inverse
of $K_+$. Since $K_+\#(-\overline{K_+})$ is a slice knot, by Lemma 3.1, 
there is a 2-disk in once punctured $-CP^2$ bounded by 
$K_-\#(-\overline{K_+})$ that represents the zero 
element. By Theorem 3.3, we have
$|1-\sigma_d(K_-\#(-\overline{K_+}))|\leq 1$, 
so we have 
$|1-\sigma_d(K_-)+\sigma_d(K_+)|\leq 1$. 
This completes the proof. $\Box$

\medskip
{\bf Proposition 3.5.} {\em Let $T(p,q)$ $(0<p<q)$ be a torus knot. 
If $T(p,q)$ is neither a trivial knot nor $T(2,3)$, then 
$\sigma_d(T(p,q))\leq -4$ for any prime integer $d$.}

\medskip
{\bf Proof.}
In \cite{PT}, J.H. Przytycki and K. Taniyama showed that, 
except for connected sums of {pretzel knots} $P(p_1,p_2,p_3)$ 
($p_1p_2p_3$ is odd), a {\em positive knot} 
can be deformed into $T(2,5)$ by changing some positive crossings 
to be negative, where a positive knot is 
a knot that has a diagram with all crossings positive. 
Since $T(p,q)$ is a prime, positive knot and 
${\rm genus}(T(p,q))\neq {\rm genus}(P(p_1,p_2,p_3))=1$, 
$T(2,5)$ is obtained from 
$T(p,q)$ by changing some positive crossings. 
Since $\sigma_d(T(2,5))=-4$ for any prime integer $d$ 
(\cite[Lemma 3.5]{Tri}), 
by Lemma 3.4, we have the conclusion. $\Box$

\medskip
{\bf Proof of Theorem 1.1.}
Note that $p\geq5$ since $T(p,q)$ is non-exceptional. 
In \cite{MM}, K. Miyazaki and K. Motegi showed that if a non-exceptional 
torus knot is obtained from a trivial knot by a single 
$(n,\omega)$-twisting, then  $|n|=1$. 
We may assume that $T(p,q)$ is obtained from a trivial knot by 
a single $(\varepsilon,\omega)$-twisting, where $|\varepsilon|=1$. 
By Lemma 3.1, there is a 2-disk $\Delta$ in a punctured 
$-\varepsilon CP^2$, $M$, such that $(\partial M,\partial\Delta)
\cong(S^3,T(p,q))$ and $\Delta$ represents 
$\omega\gamma\in H_2(M,\partial M;{\Bbb Z})$. 
If $\omega$ is divisible by a prime integer $d$, by Theorem 3.3, 
\[ 
\left|-\frac{2\varepsilon[d/2](d-[d/2])}{d^2}\omega^2+\varepsilon-
\sigma_d(T(p,q))\right|\leq 1.\]
By Propositon 3.5, $\sigma_d(T(p,q))\leq -4$. 
This gives condition (i), i.e., $\varepsilon=1$. 
So we have 
\[ 
\left|-\frac{2[d/2](d-[d/2])}{d^2}\omega^2+1-
\sigma_d(T(p,q))\right|\leq 1.\]
This implies
\[
-\frac{\sigma_d(T(p,q))}{2}\leq
\frac{[d/2](d-[d/2])}{d^2}\omega^2\leq
1-\frac{\sigma_d(T(p,q))}{2}.\]
Since $\omega$ is divisible by $d$, 
$[d/2](d-[d/2])\omega^2/d^2$ is an integer. 
This and the fact that $\sigma_d(T(p,q))$ is even (\cite[Lemma 2.16]{Tri})
give condition (iv). 
It is known that $T(p,q)$ bounds an orientable surface in $S^3$ with genus
$(p-1)(q-1)/2$. Therefore we have a closed, orientable surface with genus
$(p-1)(q-1)/2$ in $-CP^2$ that represents $\omega\gamma$. 
Thom Conjecture, which is solved by P. Kronheimer and T. Mrowka \cite{KM}, 
implies
\[\frac{(\omega-1)(\omega-2)}{2}\leq
\frac{(p-1)(q-1)}{2}.\]
Since $p-1<q-2$, we have (ii) $\omega<q$. 
Suppose that $\omega$ is even. 
By condition (iv), 
$\omega^2\geq -2\sigma(T(p,q))$. 
Since $q-1\geq p+1$, by Corollary 2.2, 
\[\omega^2\geq\left\{
\begin{array}{ll}
(p-1)(p+1)>(p-1)^2 & \mbox{ if $p$ is odd},\\
p(p+1)>p^2 & \mbox{ if $p$ is even}.
\end{array}
\right.\]
Thus we have 
$\omega> p-1$ if $p$ is odd, and 
$\omega> p$ if $p$ is even. 
Since $\omega$ is even, we have condition (iii).  $\Box$

\medskip
To prove Theorems 1.3, 1.5 and Proposition 1.8, we need the following 
theorem.

\medskip
{\bf Theorem 3.6.} (K. Kikuchi \cite{Kik})
{\em Let $M$ be a closed, oriented simply connected 
$4$-manifold with $b_2^+(M)\leq3$ and $b_2^-(M)\leq 3$. 
Let $\xi$ be a characteristic element of $H_2(M;{\Bbb Z})$. 
If $\xi$ is represented by a $2$-sphere, then 
$\xi\cdot\xi=\sigma(M)$, where $b_2^+(M)$ $($resp. $b_2^-(M))$ is 
the rank of positive $($resp. negative$)$ part of the intersection 
form of $M$.  $\Box$}

\medskip
{\bf Proof of Theorem 1.5.} (1) If $T(p,p+r)$ is obtained from a trivial knot 
by an $(n,\omega)$-twisting, then by Theorem 1.1, $n=1$ and 
$\omega>p$ if $\omega$ is even. 

Suppose $\omega$ is even. Then by Theorem 1.1 and Proposition 2.4, 
$(p-1)(p+r+1)=\omega^2,\omega^2-4$. 
Hence we have 
$(p+r/2+\omega)(p+r/2-\omega)=(r/2+1)^2$ or $(r/2+1)^2-4$. 
Since $r\geq 4$, 
\[0<\left(\frac{r}{2}+1\right)^2-4\leq\left(p+\frac{r}{2}+\omega\right)
\left(p+\frac{r}{2}-\omega\right)\leq\left(\frac{r}{2}+1\right)^2.\]
This is absurd because $p+r/2+\omega>2p+r/2 
\geq (r/2+1)^2-r/2+r/2=(r/2+1)^2$. 

Suppose $\omega$ is odd. 
Set $p=2nr+1$ $(n\geq 1)$. 
Then we have 
$O\stackrel{(1,\omega)}{\rightarrow}T(p,p+r)
\stackrel{(-1,p)}{\rightarrow}T(p,r)\cong T(r, 2nr+1)
\stackrel{(-2n,r)}{\rightarrow}T(r,1)\cong O$. 
By Lemmas 3.1 and 3.2, there is a 2-sphere in $-CP^2\#CP^2\#S^2\times S^2$ 
that represents a characteristic element 
$\omega\overline{\gamma}+p{\gamma}+r\alpha+nr\beta$. 
By Theorem 3.6, 
$-\omega^2+p^2+2nr^2=0$. 
(Note that $\omega^2=p^2+2nr^2>p^2$.)
Hence we have $(p+r/2+\omega)(p+r/2-\omega)=(r/2+1)^2-1$. 
This is absurd because $r\geq 4$, $\omega> p$  
and $2p\geq(r/2+1)^2-r/2$

(2) By the argument similar to above, we have the conclusion. 
$\Box$

\medskip
{\bf Proof of Theorem 1.3.} 
The case that $p\equiv 1$ (mod 8) is a special case of 
Theorem 1.5; see Example 1.6. Suppose that $p\equiv 3$ (mod 8). 
Set $p=8n+3$ $(n\geq 1)$. Then we have 
$O\stackrel{(1,\omega)}{\rightarrow}T(p,p+4)
\stackrel{(-1,p)}{\rightarrow}T(p,4)\cong T(4,8n+3)
\stackrel{(-2n,4)}{\rightarrow}T(4,3)\cong T(3,4)
\stackrel{(-1,3)}{\rightarrow}T(3,1)\cong O$. 
By the argument similar to that in the proof of Theorem 1.5, 
we have the conclusion. 
(Here we use Proposition 2.5 instead of Proposition 2.4.) $\Box$

\medskip
{\bf Proof of Proposition 1.8.}
Since $\sigma(T(5,8))=-20$ and 
$T(5,8)\stackrel{(-1,5)}{\rightarrow}T(5,3)
\stackrel{(-1,5)}{\rightarrow}T(5,-2)\cong T(2,-5)
\stackrel{(2,2)}{\rightarrow}T(2,-1)\cong O$. 
By the argument similar to that in the proof of Theorem 1.5, we have 
the conclusion. $\Box$

\end{document}